 \documentclass[11pt]{article}

 \usepackage{amssymb}
 \usepackage[dvips]{color}

 \newtheorem{lemma}{Lemma}[section]
 \newtheorem{theorem}{Theorem}[section]
 \newtheorem{proposition}{Proposition}[section]
 \newtheorem{corollary}{Corollary}[section]

 \def\beqlb{\begin{eqnarray}}\def\eeqlb{\end{eqnarray}}
 \def\beqnn{\begin{eqnarray*}}\def\eeqnn{\end{eqnarray*}}

 \def\proof{\noindent{\bf Proof.~~}}\def\qed{\hfill$\Box$\medskip}

 \def\<{\langle}\def\>{\rangle}

 \def\IR{{I\!\! R}}

 \def\qed{\quad$\Box$\medskip}

\def\<{\langle}\def\>{\rangle}
\def\IR{{I\!\! R}}

\def\var{\mbox{\rm var}}
 
\begin{document}
\bibliographystyle{alpha}
\noindent\underline{ December 17, 2005. Revised April 13, 2006.}

\vskip1cm

\noindent{\LARGE\bf
 A quenched CLT  for
super-Brownian
 }

\medskip

\noindent{\LARGE\bf  motion with  random immigration}

\bigskip\bigskip

\noindent{Wenming Hong\footnote{School of Mathematical Science,
Beijing Normal University, Beijing, 100875, P. R. China. 
Supported by Program for New Century  Excellent Talents in
 University (NCET) and 
  NSFC (Grant No.\,10121101) . Email: wmhong@bnu.edu.cn. }
and Ofer Zeitouni\footnote{ Department of Mathematics, University of
Minnesota and Depts. of EE and Math., Technion. This work
was partially supported by NSF grant DMS-0503775.   } }

\bigskip\bigskip

{\narrower{\narrower\small

\noindent{\bf Abstract.}  A quenched central limit theorem
is derived for the super-Brownian motion with
super-Brownian immigration, in dimension
$d\geq 4$. At the critical dimension $d=4$, the quenched and annealed
fluctuations are of the same order but are not equal.

\medskip

\noindent{\bf Key words:} Super-Brownian motion, quenched central limit
theorem,  random
immigration .

\medskip

\noindent{\bf Mathematics Subject Classification (2000):} Primary
60J80; Secondary 60F05.

\par}\par}

\bigskip
\section{Introduction and statement  of  results}

\setcounter{equation}{0}

Super-Brownian motion with super-Brownian immigration (SBMSBI, for
short), is a  superprocess in random environment, where the
environment is determined 
by an immigration process which is controlled by
the trajectory of another super-Brownian motion. Many interesting
limit properties for SBMSBI
were described under the {\it annealed} probability
(\cite{H02},  \cite{H03}, \cite{HL99} and \cite{Z05}).
In this paper, we study the central limit theorem (CLT) under
the
{\it quenched }
probability, that is, conditioned
upon a realization of the immigration process, for $d\geq 4$.

 \
To state our results and explain our motivation, we begin by
recalling the SBMSBI model (we refer to \cite{D93}  and
\cite{perkins} for a general introduction to the theory of
superprocesses). Let $C(\IR^d)$ denote the space of continuous
bounded functions on $\IR^d$. We fix a constant $p>d$ and let
$\phi_p(x):= (1+|x|^2)^{-p/2}$ for $x\in\IR^d $. Let
$C_p(\IR^d):=\{f\in C (\IR^d):\sup|f(x)|/\phi_p(x)<\infty\}$. Let
$M_p(\IR^d)$ be the space of Radon measures $\mu$ on $\IR^d$ such
that $\<\mu, f\> :=\int f(x)\mu (dx) <\infty $ for all $f\in
C_p(\IR^d)$. We endow $M_p(\IR^d)$ with the $p$-vague topology, that
is, $\mu_k\to\mu$ if and only if $\<\mu_k, f\> \to \<\mu,f\>$ for
all $f\in C_p(\IR^d)$. Then $M_p(\IR^d)$ is metrizable (\cite{I86}). 
We denote by $\lambda$  the
Lebesgue measure on $\IR^d$, and note
 that $\lambda\in M_p(\IR^d)$.

Let $S_{s,t}$ denote the heat semigroup
in $\IR^d$, that is, for $t>s$ and $f\in C(\IR^d)$,
$$S_{s,t} f(x) =\frac{1}{(2\pi (t-s))^{d/2}}
\int_{\IR^d} e^{-|y-x|^2/2(t-s)} f(y) dy\,.$$
We write $S_t:=S_{0,t}$
 and
$G$ for the corresponding potential operator, that is
$Gf=\int_0^{\infty}S_tfdt$, omitting the space variable $x$ from
the notation when no confusion may
 occur. Given $\mu\in M_p(\IR^d)$, 
a {\it super-Brownian motion} $\varrho =
(\varrho_t, P_{\mu})$ is an $M_p(\IR^d)$-valued Markov process with
$\varrho_0=\mu$ and Laplace transform given by
\begin{equation}
\label{eq-1.1}
E_{\mu}\exp\{-\<\varrho_t,f\>\} =\exp\{-\<\mu, v(t,\cdot)\>\},
\quad f\in C_p^+(\IR^d),
\end{equation}
where $v(\cdot,\cdot)$ is the unique mild solution of the
evolution equation
\begin{equation}
\label{eq-1.2}
\left\{
\begin{array}{l}
\dot{v}(t) = \frac12 \Delta v(t) - v^2(t) \\
v(0) = f\,,
\end{array}
\right.
\end{equation}
and $E_\mu$ denotes expectation with respect to $P_\mu$.

Given a super-Brownian motion $\varrho = (\varrho_t, P_{\mu})$ as
the ``environment", we will consider another super-Brownian motion
with the immigration rate controlled by the trajectory of $\varrho$,
the 
(SBMSBI) $X^{\varrho} = (X^{\varrho}_t,
P^{\varrho}_{\nu})$ with $X^{\varrho}_0=\nu$, which is again
an $M_p(\IR^d)$-valued Markov process whose {\it quenched}
probability law    is determined by
\begin{equation}
\label{eq-1.3}
E^{\varrho}_{\nu}\exp\{-\< X^{\varrho}_t,f\>\} =\exp\{-\<\nu, v(t,
\cdot)\> - \int_0^t\<\varrho_s, v(t-s, \cdot)\>ds \}.
\end{equation}
Again, $E_\nu^{\varrho}$ denotes expectations with respect to
$P_\nu^{\varrho}$.

In the following we take $\mu=\nu=\lambda$, and write $P^{\varrho}$
(resp. $P$) for $P^{\varrho}_{\lambda}$ (resp. $P_{\lambda}$). We
also use $E^{\varrho}$ and $E$ for the corresponding expectations.
This
model was considered in  \cite{HL99} and
\cite{H02,H03}, see also \cite{DGL02}, where some
interesting and new phenomena were revealed under the {\it annealed}
probability law:
$$
\mathbb{P}(\cdot):=\int P^{\varrho}(\cdot) P(d\varrho)
$$
with expectation denoted by $\mathbb{E}$. 

Our motivation for the present study is the annealed CLT derived
in \cite{HL99}, which is summarized in Theorem \ref{theo-annealed}
below.
\begin{theorem}[Hong-Li]
\label{theo-annealed}
Set
$$
\bar{a}_d(T)
 = \left\{
\begin{array}{ll}
T^{3/4}, & d=3,\\
{T^{1/2}}, & { d\geq 4, }
\end{array}
\right.
$$
and with
 $f\in C^{+}_p(R^d)$,  define 
$$
\bar Z^{\varrho_{\cdot}}_T(f):= \bar a_d(T)^{-1}\left\{ \<X^{\varrho}_T,f\> -
\mathbb{E}\<X^{\varrho}_T,f\>\right\}.
$$
Then,
$\bar Z^{\varrho_{\cdot}}_T(f)\Rightarrow \bar Z_\infty (f)$
   in distribution under
the law $\mathbb{P}$ as $T\to \infty$ ,  where $\bar Z_\infty(f)$ is a
zero mean Gaussian random variable of variance
\begin{eqnarray*}
\var(\bar Z_\infty(f))
 = \left\{
\begin{array}{ll}
\<\lambda,f\>^2/6\pi^{3/2},
 &d=3,  \\
\<\lambda,f\>^2/8\pi^2 +\<\lambda,fGf\>/2,
 &d=4,  \\
\<\lambda,fGf\>/2, &d\ge5.
\end{array}
\right.
\end{eqnarray*}
\end{theorem}
In particular, contrasting
with the standard super Brownian motion
(\cite[Theorem 5.5 and Remark 6.1]{I86}), the
SBMSBI exhibits  smoothing of the critical dimension $d=4$, since a
logarithmic term is missing in the description of the long time
behavior.

In the study of motion in
random media, differences exist between quenched and annealed
CLT behavior, and this difference is often tied to
dimension and vanishes for dimension above some critical
value. See  \cite{RS05} and \cite{Z04}  for several
examples. It is thus of interest to identify whether  similar
behavior occurs in the case of SBMSBI. Our main result, Theorem
\ref{clt-quenched} below, shows
that this is indeed the case.

Define the centered functional $Z^{\varrho_{\cdot}}_T(f)$ by
\begin{equation}
\label{eq-1.4}
Z^{\varrho_{\cdot}}_T(f):= a_d(T)^{-1}\{ \<X^{\varrho}_T,f\> -
E^{\varrho}\<X^{\varrho}_T,f\>\},
\end{equation}
where
\begin{equation}
\label{eq-1.5}
a_d(T)
 =  T^{ 1/2}, \ \ \  d\geq 4. 
 \end{equation} 
The main result of this paper is the following.
\begin{theorem}[Quenched CLT]
\label{clt-quenched}
Assume $d\geq 4$ and  $ f\in
{C^{+}_p}(\IR^d)$. Then, for $P \ a.e. \ \varrho$,
$Z^{\varrho_{\cdot}}_T(f)\Rightarrow \xi(f)$   in distribution under
the law $P^{\varrho}$ as $T\to \infty$,  where $\xi(f)$ is a
centered Gaussian variable with variance
$$\var (\xi(f))=\<\lambda, fGf\>/2.
$$
\end{theorem}
\noindent{\bf Remarks}
\begin{enumerate}
\item
As noted above,
for standard SBM in the critical dimension $d=4$, it follows from
\cite[Remark 6.1]{I86} that the  occupation measure  CLT norming is 
$(T\log
T)^{1/2}$. 
\item In \cite{H05}, the fluctuation
$b_d(T)^{-1}(E^{\varrho}\<X^{\varrho}_T,f\>-
\mathbb{E}\<X^{\varrho}_T,f\>) $ between the quenched
and annealed  means
is considered.
It is shown there that the choice
$$
b_d(T)
 = \left\{
\begin{array}{ll}
T^{(6-d)/4}  , & 3\leq d \leq 5,  \\
(\log T)^{1/2} ,
   &  d=6,\\
1,  & d\geq 7,
\end{array}
\right.
$$
leads to non-degenerate fluctuations.
Comparing Theorems \ref{theo-annealed} and \ref{clt-quenched},
one sees that in dimension $d=4$, the annealed fluctuations
consist of quenched fluctuations (around the quenched mean) and
of fluctuations of the quenched mean, and both contribute
to the annealed variance. This is not the case for $d\geq 5$:
the fluctuations of the quenched mean are of lower order and wash out
in the annealed CLT.
\item For $d=3$, an easy adaptation of our methods shows
that the statement of Theorem \ref{clt-quenched} remains true
with an almost sure statement being replaced by a
statement in probability, that is
$P^{\varrho}(
Z^{\varrho_{\cdot}}_T(f)>x)$ converges in probability, as $T\to\infty$,
to
$P(\xi(f)>x)$ for all $x$. Combined with the results
in \cite{H05}, one concludes that for $d=3$,
the quenched fluctuations
around the quenched mean are of lower order than the fluctuations
of the quenched mean. Together with Theorems \ref{theo-annealed}
and \ref{clt-quenched}, this
gives a fairly complete description
of the CLT in all dimensions $d\geq 3$.
\item A functional version of Theorem \ref{clt-quenched}
can be derived by using similar ideas. We prefer to 
bring here the shorter proof for the standard CLT.
\item Large deviations for this and related processes were studied in
\cite{IL93}, \cite{Le93} and \cite{LR95}.
\end{enumerate}

\section{Proof of Theorem \ref{clt-quenched}}
\setcounter{equation}{0}
Set $d\geq 4$ and  $f_t:=a_d(t)^{-1}f$
 with $f\in C_p^+(\IR^d)$.
For each fixed $t$,
the   mild form $v_t(r,x)$ of
equation (\ref{eq-1.2}) with $v_t(0,x)=f_t(x)$ is
\begin{equation}
\label{eq-2.1}
 v_t(r,x) =S_rf_t(x) - \int_0^rS_{r-h}v_t(h,\cdot)^2(x)dh, \quad
 0\leq r\leq t\,.
 \end{equation}
From equations
(\ref{eq-1.3}) and (\ref{eq-2.1}), it follows that
$$
E^{\varrho}\<X^{\varrho}_t,f_t\>\ =\<\lambda, S_tf_t\> +
\int_0^t\<\varrho_s, S_{t-s}f_t\>ds\,.
$$
Combined with (\ref{eq-1.4}), we get
\begin{equation}
\label{eq-2.2}
E^{\varrho} \exp\{-  Z^{\varrho}_t(f) \} =\exp\{ \<\lambda, \int_0^t
S_s v_t^2(t-s, \cdot)ds\> + \int_0^t\<\varrho_s, \int_0^{t-s} S_r
v^2_t(t-s-r, \cdot) dr\>ds \}.
\end{equation}

The proof of Theorem \ref{clt-quenched} builds upon the following two
propositions. The proof of  Proposition \ref{claim-2} will take
up most of our effort.
\begin{proposition}
\label{claim-1}
With the above notation,
$$\<\lambda, \int_0^t S_s
v_t^2(t-s, \cdot)ds\> \longrightarrow_{t\to \infty}  0\,.$$
\end{proposition}
Set $g_t(u, x)= \int_0^u S_r v_t^2(u-r, \cdot)(x) dr=
\int_0^uS_{u-r}v_t^2(r,\cdot)(x) dr$.
\begin{proposition}
\label{claim-2}
For $P$-a.e. $\varrho$,
$$\Gamma(t):=\int_0^t\<\varrho_s, g_t(t-s)\>ds \longrightarrow_{t\to\infty}
\<\lambda, fGf\>/2\,.$$
\end{proposition}

\noindent {\bf Proof of Theorem \ref{clt-quenched}}
The theorem is an immediate consequence of
(\ref{eq-2.2}), Proposition \ref{claim-1} and
Proposition \ref{claim-2}. \qed

\noindent
{\bf Proof of Proposition \ref{claim-1}}
A direct computation shows that, for any $d\geq 3$,
\begin{equation}
\label{eq-021205a}
\int_0^{\infty} ds
\< \lambda,  f S_{2s}f \> <\infty\,. \end{equation}
From (\ref{eq-2.1}), it follows that
$$\<\lambda, \int_0^t S_s v_t^2(t-s, \cdot)ds\> 
=\<\lambda, \int_0^t  v_t^2(t-s, \cdot)ds\> 
 \leq
\int_0^t\<\lambda, (S_{t-s}f_t)^2\>ds
  =  t^{-1} \int_0^t ds \< \lambda,  f S_{2s}f \>\,.
$$
Using (\ref{eq-021205a}), the result follows.\qed

\noindent {\bf Proof of Proposition \ref{claim-2}} We recall from
\cite[Theorem 3.2]{I86} that for
any $C_p(R^d)^{+}$-valued continuous path $F(s)$, the Laplace
transform  of $\int_0^t\<\varrho_s, F(t-s)\>ds$ is given by
\begin{equation}
\label{eq-2.3}
E \exp \{-\theta\int_0^t\<\varrho_s, F(t-s)\>ds\}= \exp\{-\<\lambda,
u(t,\theta; \cdot)\>\},  \ \theta >0\,,
\end{equation}
where $u(s, \theta; x)$ is the nonnegative  solution of the
following mild equation
\begin{equation}
\label{eq-2.4}
u(s,\theta;x)=\theta \int_0^sS_{s-r}F(r)(x)dr -
\int_0^sS_{s-r}u^2(r,\theta)(x)dr, \ \ 0\leq s \leq t\,.
\end{equation}
(In fact, (\ref{eq-2.3}) and (\ref{eq-2.4}) hold true for
 $|\theta|<c$ for $c$ a small enough
constant, see \cite{H03}.)
Differentiating with
respect to $\theta$ in (\ref{eq-2.3}) and (\ref{eq-2.4}),
we obtain
\begin{equation}
\label{eq-2.5}
E \left[\int_0^t\<\varrho_s, F(t-s)\>ds\right] = 
\Big\<\lambda, \int_0^t S_{t-s}F(s)ds\Big\>=
\int_0^t   \<\lambda, F(s)\>ds\,
\end{equation}
where the invariance of $\lambda$ under shifts was used 
in the second equality. Similarly,
\begin{equation}
\label{eq-2.6}
\var\left[\int_0^t\<\varrho_s, F(t-s)\>ds\right] = 2\int_0^t
\Big\<\lambda, (\int_0^s S_{s-r}F(r)dr)^2\Big\>ds\,.
\end{equation}
In the sequel, we let $A$ denote a constant whose value may change
from line to line and which may depend on the dimension and
on $f$, but not on $s,t,x$, etc.
Let us recall the  useful estimate
\begin{equation}
\label{eq-2.7}
\|S_s f\|\le A \cdot (1\wedge s^{-d/2})\,,
\end{equation}
where $\|\cdot\|$ denotes the supremum norm.
\begin{lemma}
\label{lem-2.1}
Let $F(s)=g_t(s)$. Then,
$$E \left[\int_0^t\<\varrho_s,
g_t(t-s)\>ds\right]\longrightarrow_{t\to\infty} \<\lambda, fGf\>/2\,.$$
\end{lemma}
\proof\
Note that
\begin{eqnarray*}
E \left[\int_0^t\<\varrho_s, g_t(t-s)\>ds\right] & = & \int_0^t ds
\int_0^s \<\lambda, v_t^2(r, \cdot)\>dr\\
 & = & \int_0^t ds
\int_0^s \<\lambda, (S_rf_t)^2\>dr - \int_0^t ds \int_0^s \<\lambda,
(S_rf_t)^2 -v_t^2(r, \cdot) \>dr.
\end{eqnarray*}
One has
$$
\int_0^t ds \int_0^s \<\lambda, (S_rf_t)^2\>dr = t^{-1}\int_0^t ds
\int_0^s \<\lambda, fS_{2s}f)\>dr\longrightarrow \<\lambda, fGf\>/2\,.
$$
On the other hand, by  (\ref{eq-2.1}),
\begin{eqnarray*}
 \int_0^t ds \int_0^s \<\lambda, (S_rf_t)^2 -v_t^2(r, \cdot) \>dr
 & \leq & 2 \int_0^t ds \int_0^s \<\lambda,  S_rf_t\cdot\int_0^r S_{r-h}v_t^2(h)dh
 \>dr\\
& \leq & 2 \int_0^t ds \int_0^s \<\lambda,  S_rf_t\cdot\int_0^r
S_{r-h}(S_hf_t)^2dh\>dr\\
& \leq & A t^{-3/2} \int_0^t ds \int_0^s \<\lambda, (S_rf)^2\>dr
\cdot\int_0^t (1\wedge h^{-d/2}) dh\\
&  = & A t^{-3/2}
\int_0^t ds \int_0^s \<\lambda, f S_{2r}f \>dr
\cdot\int_0^t (1\wedge h^{-d/2}) dh\\
&\leq &\frac{A}{\sqrt{t}}\cdot  \frac{1}{t}
\int_0^t ds \int_0^\infty \<\lambda, f S_{2r}f \>dr
 \end{eqnarray*}
which goes to 0 when $d\geq 3$ as $t\to\infty$ due to
(\ref{eq-021205a}); here, we used (\ref{eq-2.7}) at the
third inequality. Substituting
in (\ref{eq-2.5}), the lemma follows.\qed

We return to the proof of Proposition \ref{claim-2}.
In view of Lemma \ref{lem-2.1},
it is enough to prove that
$\Gamma(t)-E\Gamma(t)\to 0$, as $t\to\infty$, i.e.,
\begin{equation}
\label{eq-2.8}
  \int_0^t\<\varrho_s, g_t(t-s)\>ds -E
\left[\int_0^t\<\varrho_s, g_t(t-s)\>ds\right] \longrightarrow 0, \
P\ a.s.\varrho  \,.
\end{equation}
For any $n$ integer and
$n\leq t_1\leq t\leq n+1$, let $\delta= t-t_1$, and set
$\Delta
\Gamma (t_1,t) := \Gamma(t)-\Gamma(t_1)$.
Write $\Delta \Gamma (t_1,t))$ as the sum of four terms
\begin{equation}
\label{eq-2.9}
\Delta
\Gamma (t_1,t)=\Delta \Gamma_1 (t_1,t)+\Delta \Gamma_2
(t_1,t)+\Delta \Gamma_3 (t_1,t)+\Delta \Gamma_4 (t_1,t) ,
\end{equation}
 where
\begin{eqnarray*}
 \Delta \Gamma_1 (t_1,t)&:= &\int_{t_1}^t\left\<\varrho_r,
\int_0^{t-r} S_{t -r-h} v_t^2(h)dh \right\> dr,\\
\Delta \Gamma_2 (t_1,t)&:=& \int_0^{t_1}\left\<\varrho_r,
\int_{t_1-r}^{t-r} S_{t -r-h} v_t^2(h)dh \right\> dr,\\
 \Delta
\Gamma_3 (t_1,t) &:=&\int_0^{t_1}\left\<\varrho_r, \int_0^{t_1-r}
 S_{t-r-h}[v_{t }^2(h)-v_{t_1}^2(h)]dh \right\> dr ,\\
\Delta \Gamma_4 (t_1,t)&:= & \int_0^{t_1}\left\<\varrho_r,
\int_0^{t_1-r} [S_{t-r-h}-S_{t_1-r-h}]v_{t_1}^2(h)dh \right\> dr.\\
 \end{eqnarray*}
We estimate separately the
moments of centered versions of $\Delta \Gamma_i  (t_1,t) $.
\begin{lemma}
\label{lem-step1}
$$\var[\Delta \Gamma_1 (t_1,t)]  \leq A
\delta^2 n^{-2} \,.$$
\end{lemma}
\proof\ Recall that $\Delta \Gamma_1 (t_1,t)=
\int_{t_1}^t\left\<\varrho_r, g_t(t-r) \right\> dr$. We have, again
from  
\cite[Theorem 3.2]{I86}, that for
$\theta\geq 0$,
\begin{eqnarray*}
E \exp \left\{-\theta \int_{t_1}^t\left\<\varrho_r, g_t(t-r)
\right\> dr \right\}&=&E \exp\left\{- \left\<\varrho_{t_1}, u(t_1,
t,\theta;\cdot) \right\> \right\}\\
&=& \exp\left\{- \left\<\lambda, w(0,t_1, \theta;\cdot) \right\>\,,
\right\}
\end{eqnarray*}
where $u(s, t,\theta;\cdot)$ is the nonnegative  solution of the
following mild equation
$$
u(s, t,\theta;x)=\theta \int_s^tS_{s,r}g_t(t-r)(x)dr -
\int_s^tS_{s,r}u^2(r, t,\theta)(x)dr, \ \ t_1\leq s \leq t\,,
$$
and $w(s, t_1,\theta;x)$  is the nonnegative  solution of the
following mild equation
$$
w(s, t_1,\theta;x)=  S_{s,t_1}u(t_1, t,\theta;\cdot)(x) -
\int_s^{t_1}S_{s,r}w^2(r, t_1,\theta)(x)dr, \ \ 0\leq s \leq t_1\,.
$$
Obviously,
$$\var[\Delta \Gamma_1 (t_1,t)] = - \frac{\partial^2 \<\lambda,
w(0,t_1,\theta;\cdot)\>}{\partial \theta^2}\mid_{\theta=0}\,.$$
Performing the differentiation and using that
$u\mid_{\theta=0}=w\mid_{\theta=0}=0$,
we obtain
\begin{eqnarray*}
& &\var[\Delta \Gamma_1 (t_1,t)] \\
&=& 2 \Big\<\lambda, \int_{t_1}^tS_{0 , s}\Big[ \int^t_sS_{s,r}
g_t(t-r)dr\Big]^2ds\Big\> + 2\Big\<\lambda,
\int_0^{t_1}S_{0,s}\Big[\int_{t_1}^tS_{s,r}g_t(t-r)dr\Big]^2ds\Big\>\\
&\leq & 2t^{-2}\left[ \Big\<\lambda, \int_{t_1}^tS_{s   }\Big[ \int^t_s
\int_0^{t-r}S_{t-s-l}(S_lf)^2dldr\Big]^2ds\Big\> +  \Big\<\lambda,
\int_0^{t_1}S_{ s}\Big[\int_{t_1}^t \int_0^{t-r}S_{t-s-l}(S_lf)^2dldr\Big]^2ds
\Big\>\right]\\
&\leq & At^{-2}\left[ \Big\<\lambda, \int_{t_1}^t [ (t-s)S_{t-s}f]^2ds\Big\>
+  \Big\<\lambda,
\int_0^{t_1} \Big[\int_{t_1}^t  S_{t-s}fdr\Big]^2ds\Big\>\right]\\
&\leq & An^{-2}\delta^2,
\end{eqnarray*}
in which (\ref{eq-2.7}) has been used several times. \qed

\begin{lemma}
\label{lem-step2}
With the above notation,
$$ E[\Delta \Gamma_2 (t_1,t)-E\Delta \Gamma_2
(t_1,t)]^4  \leq A \delta^2 n^{-2} \,.$$
\end{lemma}
\proof\
We have
$\Delta \Gamma_2 (t_1,t)= \int_0^{t_1}\left\<\varrho_r,F_t(t_1-r)
\right\> dr$ where
 $F_t(r):=\int_r^{\delta+r}S_{\delta+r-l}v_t^2(l)dl$, then
$\Delta \Gamma_2 (t_1,t)= \int_0^{t_1}\left\<\varrho_r,F(t_1-r)
\right\> dr$. Let
$$u^{(i)}(r,x):=\frac{\partial^i u(r,  x,
\theta)}{\partial \theta^i}|_{\theta=0}\,, i=1,2,3\,.$$
Differentiating with respect to $\theta$ in (\ref{eq-2.3}) and (\ref{eq-2.4}),
and using again that $u\mid_{\theta=0}=0$, we obtain
\begin{eqnarray}
\label{eq-021205b}
 & & E[\Delta \Gamma_2 (t_1,t)-E\Delta \Gamma_2 (t_1,t)]^4\nonumber \\
 &=&3\left(\int_0^{t_1}\<\lambda, u^{(1)}(r)^2\>dr\right)^2+3\int_0^{t_1}
\<\lambda,
 u^{(2)}(r)^2\>dr + 4\int_0^{t_1}\<\lambda,
 u^{(1)}(r)u^{(3)}(r) \>dr\nonumber \\
 &:=&3 I^2+3J+4K,
\end{eqnarray}
where for $0\leq r\leq t_1$,
\begin{eqnarray*}
u^{(1)}(r,x)&=&\int_0^{r}S_{r-s}F_t(s)ds =
\int_0^{r}S_{r-s}\int_s^{\delta+s}S_{\delta+s-l}v_t^2(l)dlds \\
   &\leq & A t_1^{-1}\int_0^{r}S_{r-s}\int_s^{\delta+s}S_{\delta+s-l}(S_lf)^2dlds\\
   &\leq &  A t_1^{-1}\delta \cdot rS_{\delta+r}f\\
|u^{(2)}(r,x)|&=&\Big|-2
\int_0^{r}S_{r-s}u'(s)^2ds\Big|\leq A t_1^{-2}
\delta^2\cdot \int_0^{r}S_{r-s}(sS_{\delta+s}f)^2ds\\
         & \leq & A t_1^{-2}\delta^2\cdot rS_{\delta+r}f \\
|u^{(3)}(r,x)|&=&\Big|-6\int_0^{r}S_{r-s}u'(s)u''(s)ds\Big|
\leq A t_1^{-3}\delta^3\cdot \int_0^{r}S_{r-s}(sS_{\delta+s}f)^2ds\\
         & \leq & A t_1^{-3}\delta^3\cdot rS_{\delta+r}f\,. \\
\end{eqnarray*}
Thus, we obtain, using that  $d\geq 4$,
\begin{eqnarray*}
I&\leq & A  t_1^{-2}\delta^2\cdot
\int_0^{t_1}\<\lambda,(rS_{\delta+r}f)^2\>dr \leq A
t_1^{-2}\delta^2\cdot \int_0^{t_1} r^2(1\wedge (\delta+r)^{-d/2})dr
\leq A  t_1^{-1}\delta^2.
\end{eqnarray*}
Similarly, we have
 \begin{eqnarray*}
 J=\int_0^{t_1}\<\lambda, u^{(2)}(r)^2\>dr
 &\leq & A  t_1^{-4}\delta^4\cdot
\int_0^{t_1}\<\lambda,(rS_{\delta+r}f)^2\>dr \leq A
t_1^{-3}\delta^4,
\end{eqnarray*}
and
\begin{eqnarray*}
 K=\int_0^{t_1}\<\lambda,
 u^{(1)}(r)u^{(3)}(r) \>dr
 &\leq & A  t_1^{-4}\delta^4\cdot
\int_0^{t_1}\<\lambda,(rS_{\delta+r}f)^2\>dr \leq A
t_1^{-3}\delta^4.
\end{eqnarray*}
Substituting in
(\ref{eq-021205b}) completes the proof. \qed

\begin{lemma}
\label{lem-step3}
With the above notation,
$$\var[\Delta \Gamma_3  (t_1,t)]  \leq A
\delta^2 n^{-3} \,.$$
\end{lemma}

\proof\
We begin by considering
the difference $v_t(r,x)-v_{t_1}(r,x)$.
 From
(\ref{eq-2.1}), we have
\begin{equation}
\label{eq-2.10}
v_t(r,x)-v_{t_1}(r,x)= S_rf_t(x)-S_rf_{t_1}(x)-
\int_0^rS_{r-h}[v_t^2(h,\cdot)-v_{t_1}^2(h,\cdot)](x)dh.
\end{equation}
A direct computation reveals that
  $ \|S_rf_t -S_rf_{t_1}\|\leq A
\delta t_1^{-3/2}(1\wedge r^{-d/2})$. Since
 $v_t(r,x)\leq
S_rf_t$, it follows that  $\|S_{r-h}[v_t^2(h,\cdot)-v_{t_1}^2(h,\cdot)]\|\leq
\|2[v_t (h,\cdot)-v_{t_1} (h,\cdot)]S_r f_{t_1}\|\leq A
t_1^{-1/2}(1\wedge r^{-d/2})\|v_t (h,\cdot)-v_{t_1} (h,\cdot)\|.$
Thus from (\ref{eq-2.10}) we get,
$$
\|v_t(r,\cdot)-v_{t_1}(r,\cdot)\|\leq A \delta t_1^{-3/2}(1\wedge
r^{-d/2})+ A t_1^{-1/2}(1\wedge r^{-d/2})\int_0^r\|v_t
(h,\cdot)-v_{t_1} (h,\cdot)\|dh\,.
$$
Writing $a_r=\|v_t(r,\cdot)-v_{t_1}(r,\cdot)\|\geq 0$,
$b_r=A \delta t_1^{-3/2}(1\wedge
r^{-d/2})\geq 0$ and $c_r=A t_1^{-1/2}(1\wedge r^{-d/2})\geq 0$, we thus have
$$ a_r\leq b_r+c_r\int_0^r a_s ds\,.$$
By a version of Gronwall's inequality,
$$ a_r\leq b_r+c_r\int_0^r e^{\int_s^t c_u du} b_s ds\,.$$
(This can be seen by setting $z_r=\int_0^r a_s ds$
and noting that $z_r$ satisfies the differential inequality
$d z_r/dr \leq b_r+c_r z_r$, with $z_0=0$.)
Thus,
\begin{eqnarray*}
& & \|v_t(r,\cdot)-v_{t_1}(r,\cdot)\|\\
 &\leq  & A \delta
t_1^{-3/2}(1\wedge r^{-d/2})+ A^2 t_1^{-1/2}(1\wedge r^{-d/2})\int_0^r
\delta t_1^{-3/2}(1\wedge
s^{-d/2})\exp\{At_1^{-1/2}\int_s^r (1\wedge u^{-d/2})du\}ds\\
& \leq &    A \delta t_1^{-1}(1\wedge r^{-d/2}).
\end{eqnarray*}
Once more by (\ref{eq-2.10}) we have
\begin{eqnarray*}
& & |v_t(r,x)-v_{t_1}(r,x)|\\
 &\leq & |S_rf_t(x)-S_rf_{t_1}(x)|+
\int_0^r||v_t (h,\cdot)-v_{t_1} (h,\cdot)||\cdot |S_{r-h}[v_t
(h,\cdot)+v_{t_1} (h,\cdot)] (x)|dh\\
&\leq & \delta t_1^{-3/2}S_rf(x)+A \delta t_1^{-3/2}S_rf(x) \cdot \int_0^r
(1\wedge h^{-d/2})dh\\
&\leq &  A \delta t_1^{-3/2}S_rf(x).
\end{eqnarray*}
Now we can estimate the variance of $\Gamma_3
(t_1,t)$. By (\ref{eq-2.6}) with $F(r)=\int_0^r S_{t-t_1+r-l}[v_{t
}^2(l)-v_{t_1 }^2(l)]dl$,
\begin{eqnarray*}
\var[\Delta \Gamma_3   (t_1,t)]  &=&2\int_0^{t_1}\left\<\lambda,
\left[\int_0^{ r}S_{r-h}   \int_0^h S_{t-t_1+h-l}[v_{t }^2(l)-v_{t_1
}^2(l)]dldh\right]^2 \right\> dr\\
&\leq & A \int_0^{t_1}\left\<\lambda, \left[\int_0^{ r} \int_0^h
S_{\delta+r-l}|v_{t } (l)-v_{t_1 } (l)||v_{t } (l)+v_{t_1
} (l)|dldh\right]^2 \right\> dr\\
&\leq & A \delta^2 t_1^{-4}\int_0^{t_1}r^2[1\wedge (r+\delta)^{-d/2}]dr\\
&\leq & A \delta^2 n^{-3}.
\end{eqnarray*}
This completes the proof of the lemma.
 \qed

Before providing an estimate on the moments of $\Delta \Gamma_4(t_1,t)$,
we need an a-priori simple estimate on time differences of
the heat kernel $p(t,x,y)=(2\pi t)^{-d/2} \exp(-|x-y|^2/2t)$.
Since we did not find a direct reference for it,
we provide the proof.
\begin{lemma}
\label{lem-2.2}
There is a constant $A$ such  that for any $t\geq \tau>0$,   we have
\begin{equation}
\label{eq-2.11}
\sup_{0<s\leq \tau\leq t} s^{-1}\left|p(t+s,x,y)-p(t
,x,y)\right|\leq A \tau^{-1} [p(t+2\tau,x,y)+p(t ,x,y)],
\end{equation}
\end{lemma}
\proof
Consider first $\tau=1$. Let $z=|x-y|$, two
cases should be considered:

\noindent{\it Case 1}:  $ z^2< 2d(t+1)$. Note that
$$
\left|p(t+s,x,y)-p(t
,x,y)\right|=p(t,x,y)\left|\exp\left\{\frac{z^2s}{2t(t+s)}\right\}
\left(\frac{t}{t+s}\right)^{d/2}-1\right|
$$
and $\exp\left\{\frac{z^2s}{2t(t+s)}\right\}=1+sR_1(s,t,z)$,
$\left(\frac{t}{t+s}\right)^{d/2}=
\left(1-\frac{s}{t+s}\right)^{d/2}=1+sR_2(s,t,z)$,
where $ R_1(s,t,z)$, $R_2(s,t,z)$ are bounded by a constant
when $ z^2< 2d(t+1)$, $0<s\leq
1\leq t$.  Thus,
we get
\begin{equation}
\label{eq-2.12}
\sup_{0<s\leq 1}\sup_{z^2< 2d(t+1)} s^{-1}\left|p(t+s,x,y)-p(t
,x,y)\right|\leq A    p(t ,x,y)\,.
\end{equation}

\noindent{\it Case 2}:  $ z^2\geq 2d(t+1)$. Since $ \frac{
\partial}{\partial
t}p(t,z)=p(t,z)\left[-\frac{d}{2t}+\frac{z^2}{2t^2} \right] $,
\begin{eqnarray*}
\left|p(t+s,z)-p(t,z)\right|
 &=&\left|\int_0^s
p(t+u,z)\left[-\frac{d}{2(t+u)}+\frac{z^2}{2(t+u)^2} \right]du\right|\\
&\leq &\int_0^s p(t+u,z) \frac{z^2}{2(t+u)^2} du\,,
 \end{eqnarray*}
where the inequality uses that
$\left|-\frac{d}{2(t+u)}+\frac{z^2}{2(t+u)^2} \right|\leq
\frac{z^2}{2(t+u)^2}$ when  $ z^2\geq 2d(t+1)$. But
\begin{eqnarray*}
\frac{z^2}{2(t+u)^2}p(t+u,z)
&=&p(t+2,z)\left(\frac{t+2}{t+u}\right)^{d/2}\frac{z^2}{(t+u)^2}
\exp\left\{-\frac{z^2}{2}\left[\frac{2-u}{(t+u)(t+2)}\right]\right\}\\
&\leq& A \cdot
p(t+2,z)\frac{z^2}{(t+u)^2}
\exp\left\{-\frac{z^2}{2}\left[\frac{2-u}{(t+u)(t+2)}\right]\right\},
 \end{eqnarray*}
(note that $0<u\leq s\leq 1\leq t $) and
$$
\sup_{z^2\geq
2d(t+1)}\frac{z^2}{(t+u)^2}
\exp\left\{-\frac{z^2}{2}\left[\frac{2-u}{(t+u)(t+2)}\right]\right\}
 <\infty.
$$
So
\begin{equation}
\label{eq-2.13}
\sup_{0<s\leq 1}\sup_{z^2\geq 2d(t+1)} s^{-1}\left|p(t+s,x,y)-p(t
,x,y)\right|\leq A    p(t+2 ,x,y). \end{equation}
Combining (\ref{eq-2.12}) and (\ref{eq-2.13}) we obtain (\ref{eq-2.11})
when $\tau=1$. For
general $\tau>0$, we use the scaling properties of
$p(t,z)$. We have
\begin{eqnarray*}
s^{-1}\left|p(t+s,z)-p(t ,z)\right|
&=&\tau^{-d/2}s^{-1}\left|p(\tau^{-1}(
t+s),\tau^{-1/2}z)-p(\tau^{-1}t ,\tau^{-1/2}z)\right|\\
&=&\tau^{-d/2}\tau^{-1}\left[\tau s^{-1}\left|p(\tau^{-1}(
t+s),\tau^{-1/2}z)-p(\tau^{-1}t ,\tau^{-1/2}z)\right|\right]\\
&\leq &A \tau^{-d/2}\tau^{-1}\left[  p(\tau^{-1}
t+2),\tau^{-1/2}z)+p(\tau^{-1}t ,\tau^{-1/2}z) \right]\\
&\leq &A  \tau^{-1}\left[  p( t+2\tau ), z)+p( t , z) \right],
\end{eqnarray*}
where the third inequality follows from the case already considered
because
$0<\tau^{-1} s\leq 1\leq\tau^{-1}t $. This complete the proof. \qed

\begin{lemma}
\label{lem-step4}
With the notation above, we have
$$ E[\Delta \Gamma_4  (t_1,t)-E\Delta \Gamma_4
(t_1,t)]^4 \leq A \delta^2 n^{-2} $$\end{lemma}
\proof
The formula for the
 fourth  moment of $\Delta \Gamma_4
(t_1,t)$
is  as in the proof of Lemma \ref{lem-step2}, except that
the function $F_t(r)$ is replaced by
the function
 $\tilde F_{t_1}(r):=\int_0^r[S_{\delta+r-l}-S_{
r-l}]v_{t_1}^2(l)dl$, and $\Delta \Gamma_4  (t_1,t)= \int_0^{t_1}
\<\lambda, \tilde F_{t_1}(r)\>dr$. Recalling that
$\delta=t-t_1$, we obtain for
$0\leq r\leq t_1$ that,
\begin{eqnarray*}
|u^{(1)}(r,x)|&=&\left|\int_0^{r}S_{r-s}\tilde F_{t_1}(s)ds \right|=
\left|\int_0^{r}S_{r-s}\int_0^s
[S_{\delta+s-l}-S_{ s-l}]v_{t_1}^2(l)dlds \right|\\
&\leq & \int_0^r l \left|[S_{\delta+ l}-S_{  l}]v_{t_1}^2(r-l)\right|dl \\
&\leq & \int_0^{\delta} l   [S_{\delta+ l}+S_{  l}](S_{(r-l)}f_{t_1})^2 dl
+ A\int_{\delta}^r \delta  [S_{3 l}+S_{  l}](S_{(r-l)}f_{t_1})^2 dl\\
&\leq & A \delta t_1^{-1}\left[\delta(S_{\delta+r}+S_{r})f
+\int_{\delta}^r   (S_{r+2 l}+S_{r})f \cdot (1\wedge (r-l)^{-d/2})dl
\right],
\end{eqnarray*}
where the third step is from Lemma \ref{lem-2.2} (with $\tau=t=l, s=\delta$
there). By a similar calculation we get
\begin{eqnarray*}
 |u^{(2)}(r,x)|&=&\left|-2\int_0^{r}S_{r-s}u^{(1)}(s)^2ds\right| \\
         & \leq & A t_1^{-2}\delta^2\cdot
         \left[\delta(S_{\delta+r}+S_{r})f
+\int_0^rds \int_{\delta}^s   (S_{r+2 l}+S_{r})f \cdot (1\wedge
(s-l)^{-d/2})dl \right],\\
|u^{(3)}(r,x)|&=&\left|-6\int_0^{r}S_{r-s}u^{(1)}(s)u^{(2)}(s)ds\right| \\
         & \leq & A t_1^{-3}\delta^3\cdot \left[\delta(2r^{3/2}+\delta)(S_{\delta+r}+S_{r})f
+r^{1/2}\int_0^rds \int_{\delta}^s   (S_{r+2 l}+S_{r})f \cdot
(1\wedge (s-l)^{-d/2})dl \right],
\end{eqnarray*}
and the estimate (\ref{eq-2.7}) was used many times. Then
\begin{eqnarray*}
I&=& \int_0^{t_1}\<\lambda, u^{(1)}(r)^2\>dr\\
&\leq & A \delta^2
t_1^{-2}\int_0^{t_1}\left\<\lambda,\left[\delta(S_{\delta+r}+S_{r})f
+\int_{\delta}^r   (S_{r+2 l}+S_{r})f \cdot (1\wedge (r-l)^{-d/2})dl
\right]^2 \right\>dr\\
&\leq   & A\delta^2 t_1^{-1},
\end{eqnarray*}
and $J =  3\int_0^{t_1}\<\lambda,
 u^{(2)}(r)^2\>dr \leq A\delta^4 t_1^{-2},$
$K  =  \int_0^{t_1}\<\lambda,
 u^{(1)}(r)u^{(3)}(r) \>dr \leq A\delta^4 t_1^{-2}.
$ So
\begin{eqnarray*}
   E[\Delta \Gamma_4  (t_1,t)-E\Delta
\Gamma_4 (t_1,t)]^4  = 3 I^2+3J+4K \leq A\delta^4
  n^{-2},
\end{eqnarray*}
which completes the proof.
\qed

\noindent
We return to the proof of Proposition \ref{claim-2}.
Let $\overline{\Gamma}(t)
:=\Gamma(t)-E\Gamma(t)$ denote the  centered $\Gamma(t)$, and
define
$\Delta
\overline{\Gamma_i}$ similarly. For any $\varepsilon>0$
and $\alpha\in (1,4/3)$,
\begin{eqnarray*}
& &P(\max_{n^{\alpha}\leq t\leq (n+1)^{\alpha}}
\left|\overline{\Gamma}(t)
-\overline{\Gamma}(n^\alpha) \right|>\varepsilon)
\\
& \leq & \sum_{k=1}^{\infty} P(\max_{0\leq j\leq 2\alpha n^{\alpha-1}
2^k}\left|\overline{\Gamma}(n^\alpha+2^{-k}(j+1))
-\overline{\Gamma}(n^\alpha+2^{-k} j ) \right|>
\frac{\epsilon}{4\alpha n^{\alpha-1} k^2})\\
& = & \sum_{k=1}^{\infty} P(\max_{0\leq j\leq 2\alpha n^{\alpha-1}
2^k }\left|\Delta\overline{\Gamma}(n^\alpha+2^{-k} j, n^\alpha+
2^{-k}(j+1))
 \right|>
\frac{\epsilon}{4\alpha n^{\alpha-1} k^2})\\
& \leq & \sum_{k=1}^{\infty}\sum_{i=1}^{4}  
2\alpha n^{\alpha-1}2^k
\max_{0\leq j\leq 2\alpha n^{\alpha-1}2^k
} P(\left| \Delta\overline{\Gamma_i} (n^\alpha+2^{-k} j,
n^\alpha+2^{-k}(j+1))
 \right|>
\frac{\epsilon}{16\alpha n^{\alpha-1} k^2})\,.
\end{eqnarray*}
By Chebyshev's inequality and Lemmas \ref{lem-step1}
and \ref{lem-step3}, for $i=1,3$,
$$P(\left| \Delta\overline{\Gamma_i} (n^\alpha+2^{-k} j, n^\alpha
+2^{-k}(j+1))
 \right|>
\frac{\epsilon}{16\alpha n^{\alpha-1} k^2}
)\leq 
A_\alpha \epsilon^{-2} k^4 2^{-2k} n^{-2}\,.$$
 Similarly,
using Lemmas \ref{lem-step2} and \ref{lem-step4}, we obtain for
 $i=2,4$,
$$P(\left| \Delta\overline{\Gamma_i} (n^\alpha
+2^{-k} j, n^\alpha
+2^{-k}(j+1))
 \right|>
\frac{\epsilon}{16\alpha n^{\alpha-1} k^2})
\leq 
A_\alpha \epsilon^{-4} k^8 2^{-2k} n^{2\alpha-4}.$$
Thus, adjusting the value of $A_\alpha$, using that
$2\alpha-4>-2$,
\begin{eqnarray*}
 P(\max_{n^{\alpha}\leq t\leq (n+1)^{\alpha}}\left|\overline{\Gamma}(t)
-\overline{\Gamma}(n^\alpha) \right|>\varepsilon)
 &  \leq & A_\alpha {\varepsilon^{-4}}n^{3\alpha-5}
\sum_{k=1}^{\infty}       { k^{
 8}}2^{- k} \leq A_\alpha {\varepsilon^{-4}}n^{3\alpha-5} \,.
\end{eqnarray*}
By the Borel-Cantelli Lemma, we get $\max_{n^{\alpha}\leq t\leq
(n+1)^{\alpha}}\left|\overline{\Gamma}(t) -\overline{\Gamma}(n^\alpha)
\right|\to 0$,  $P-a.s.$. Thus, the proposition follows once we prove
that
\begin{equation}
\label{eq-2.14}
\overline{\Gamma}(n^\alpha)\longrightarrow 0 \ \ \ P-a.s.
\end{equation}
Recall that $E[\overline{\Gamma}(n^\alpha)]=0$, and by (\ref{eq-2.6}),
\begin{eqnarray*}
 \var[\overline{\Gamma}(n^\alpha)]&=
 & \var\left[\int_0^{n^\alpha}\<\varrho_s,
g_{n^\alpha}(n^\alpha-s)\>ds\right]
  =  2\int_0^{n^\alpha}
\Big\<\lambda, (\int_0^s S_{s-r}g_{n^\alpha}(r)dr)^2\Big\>ds\\
&\leq & 2\int_0^{n^\alpha}
\Big\<\lambda, \Big[\int_0^s S_{s-r}
(\int_0^rS_{r-l}(S_lf_{n^\alpha})^2dl)dr\Big]^2
\Big\>ds\\
&\leq & A_\alpha \cdot n^{-2\alpha} \int_0^{n^\alpha}
\Big\<\lambda, [sS_sf\cdot  \int_0^{\infty}(1\wedge 
l^{-d/2}) dl]^2\Big\>ds\\
&\leq & A_\alpha \cdot n^{-2\alpha} 
\int_0^{n^\alpha} s^2 (1\wedge s^{-d/2})ds  \\
&\leq & A_\alpha \cdot n^{-\alpha}.
\end{eqnarray*}
Thus for any $\varepsilon>0$,
$$
\sum_{n=1}^{\infty} P[|\overline{\Gamma}(n^\alpha)|>\varepsilon] \leq
A_\alpha
\sum_{n=1}^{\infty} \varepsilon^{-2}n^{-\alpha}
 < \infty,$$  and (\ref{eq-2.14})
follows by the Borel-Cantelli Lemma. \qed

 \noindent {\bf Acknowledgment}  This work was done during a visit
of Wenming Hong to the University of Minnesota. He would like to
thank the Department of Mathematics, University of Minnesota for its 
hospitality during this visit.


\noindent
\begin{thebibliography}{999999}












\bibitem[D93]{D93}
Dawson, D.A.,
Measure-valued Markov  processes, In: {\it Lect. Notes. Math.}  1541,
1-260 (1993), Springer-Verlag.





\bibitem[DGL02]{DGL02}
Dawson, D.A., Gorostiza, L.G., Li, Z.H.,  { Non-local
branching superprocesses and some related models}, {\it Acta
Applicandae Mathematicae},  74 (2002), 93--112.








\bibitem[H02]{H02} Hong, W.M., Longtime behavior
for the occupation time of super-Brownian motion
with random immigration, {\it Stochastic Process. Appl.} 102 (2002),
43--62.




\bibitem[H03]{H03} Hong, W.M., Large deviations for
the super-Brownian motion
with super-Brownian immigration, {\it Journal of Theoretical
Probability}, 16 (2003), 899-922.

\bibitem[H05]{H05} Hong, W.M.,  Quenched mean limit theorems for
the super-Brownian motion with super-Brownian immigration, {\it
Infinite Dimensional Analysis, Quantum Probability and Related
Topics},  8 (2005), 383-396. 


\bibitem[HL99]{HL99} Hong, W.M. and Li, Z.H.,
 A central limit theorem for the super-Brownian
motion with super-Brownian immigration, {\it J. Appl. Probab.} 36 (1999),
1218-1224.





\bibitem[I86]{I86} Iscoe, I., 
A weighted occupation time for a class of measure-valued
critical branching Brownian motion, {\it Probab. Th. Rel. Fields}
71 (1986), 85-116.


\bibitem[IL93]{IL93} Iscoe, I., Lee, T.Y., 
 Large deviations for occupation times of measure-valued branching Brownian
motions, {\it Stoch. Stoch. Rep.}  45 (1993), 177-209.




\bibitem[L93]{Le93} Lee, T.Y., Some
 limit theorems for super-Brownian motion and semilinear differential
equations. {\it Ann. Probab.} 21 (1993), 979-995.

\bibitem[LR95]{LR95} Lee, T.Y. and Remillard, B, 
  Large deviation for three dimensional super-Brownian motion.
{\it Ann. Probab.}
 23 (1995), 1755-1771.

\bibitem[P02]{perkins} Perkins, E. A., Dawson-Watanabe superprocesses and
measure-valued diffusions, {\it Lecture Notes  Math.} 1781 (2002), 
 132--318,
Springer-Verlag.

\bibitem[RS05]{RS05}  
Rassoul-Agha, F.,   Sepp\"{a}l\"{a}inen, T.:
 An almost sure invariance principle for random
walks in a space-time random environment,
{\it   Probab. Theory Related Fields} {\bf  133}  ,
pp.  299--314 (2005).

\bibitem[Ze04]{Z04} Zeitouni, O., {Random walks
in random environment},
{ \it Lecture notes in Math.} 1837 (2004), 193--312, Springer.

\bibitem[Zh05]{Z05} Zhang, M., Functional central
limit theorems  for the super-Brownian motion with super-Brownian
immigration, {\it Journal of Theoretical Probability}, 18 (2005), 665-685.

\end{thebibliography}
\end{document}